\documentclass[preprint,12pt]{elsarticle}

\usepackage{amssymb}

\usepackage{graphicx}
\usepackage{dcolumn}
\usepackage{bm}
\usepackage{amsmath}
\usepackage{amssymb}
\usepackage{latexsym}
\usepackage{epsfig}
\usepackage{amsbsy}
\usepackage{array}
\usepackage{amssymb}
\usepackage{setspace}
\usepackage{bm}
\usepackage{amsthm}
\usepackage{subfigure}

\def\sint{\ifmmode{- \!\!\!\!\!\! \int}
    \else{\hbox{$- \!\!\!\! \int \ $}}\fi}

\newtheorem{theorem}{Theorem}
\begin{document}

\begin{frontmatter}



\title{Approximating rational B\'{e}zier curves by constrained B\'{e}zier curves of arbitrary degree}

\author[label1,label2]{Mao Shi}
\author[label2]{Jiansong Deng}
\address[label1]{College of Mathematics and Information Science of Shaanxi Normal University,Xi'an 710062, China}
\address[label2]{Department of Mathematics, University of Science and Technology of China, Hefei, Anhui 230026, China}
\cortext[cor1]{ Email: shimao@snnu.edu.cn}

\begin{abstract}
In this paper, we propose a method to obtain a constrained approximation of a rational B\'{e}zier curve by a polynomial B\'{e}zier curve. This problem is reformulated as an approximation problem between two polynomial B\'{e}zier curves based on weighted least-squares method, where weight functions $\rho(t)=\omega(t)$ and $\rho(t)=\omega(t)^{2}$ are studied respectively. The efficiency of the proposed method is tested using some examples.
\end{abstract}


\end{frontmatter}

\section{Introduction}
Rational B\'{e}zier curves play a significant role in Computer Aided Design systems. However,  the forms of  derivatives are quite complex and integral expressions may not exist for high-degree rational B\'{e}zier curves, so the problem of approximating rational functions with polynomials has been raised and studied. Sederberg and Kakimoto \citep {Sederberg} presented the hybrid polynomial approximation to rational curves for the first time in 1991. Wang et al.\citep {Wang}\citep{Wang08} presented Hermite polynomial approximations to rational B\'{e}zier curves and investigated the convergence condition for the polynomial approximation of rational functions and rational curves. Floater \citep {Floater} constructed a high-order approximation of rational curves using polynomial curves. Lee \citep {Lee} converted a polynomial approximation of a rational B\'{e}zier curve to control points approximation of two B\'{e}zier  curves. An approximate solution was obtained by  least-squares method.
Huang \citep{Huang} presented a simple method for approximating a rational B\'{e}zier curve with a B\'{e}zier curve sequence based on degree elevation. Recently,  sample-based polynomial approximation of rational B\'{e}zier curves was investigated by Lu \citep{Lu}, whose errors   smaller  than Huang's. Hu \citep{Hu13} proposed a reparameterization-based method for polynomial approximating rational B\'{e}zier curves with constraints.

In this paper, we propose a method to obtain a constrained approximation of a rational B\'{e}zier curve by a polynomial B\'{e}zier curve based on  weighted least-squares \cite{Isaacson}. The main idea is converting a polynomial approximation of a rational B\'{e}zier curve to an approximation between two  B\'{e}zier curves   so that control points of a B\'{e}zier curve can be obtained using  linear equations.  Compared with the related works in paper \citep{Lee} , \citep{Lu} and \citep{Hu13}, our method is fast, stable and  highly accurate  to obtain some higher degree  B\'{e}zier curves as well as don't need  fix parameter $\lambda$.

The  paper is structured as follows. Section 2  presents some basic concepts and properties regarding the problem of the constrained B\'{e}zier approximation  of a rational B\'{e}zier curve. Section 3 brings a complete solution
to the problem formulated above in the $L_{2}$ norm. Section 4 presents some numerical examples  to verify the accuracy and effectiveness of the method.

\section{Preliminaries}
\label{Section:theory}
\subsection{Definitions and Properties}
A standard-form  $n$th degree rational B\'{e}zier curve is defined as follows {\cite {Farin}}:
 \begin{eqnarray}
\textit{\textbf{P}}(t)=\frac{\textbf{\textit{x}}(t)}{\omega (t)}=\frac{\sum \limits_{i=0}^n \omega_{i}\textit{\textbf{p}}_{i}B_{i}^{n}(t)}{\sum \limits_{i=0}^{n}\omega_{i}B_{i}^{n}(t)},\ t\in [0, 1],
\label{eq:energy}
\end{eqnarray}
 where $B_{i}^{n}(t)={n \choose i}t^{i}(1-t)^{n-i}$ are the Bernstein polynomials, $\textit{\textbf{p}}_{i}\in \mathbb{R}^{d}, i=0,1,...,n,$ are  control points, and $\omega_{i}\in \mathbb{R}^{+},  i=0,1,...,n, \omega_{0}=\omega_{n}=1$,  are the  scalar weights. Clearly, polynomial function $\omega (t)>0 $.

 The following are pertinent theorems used in this paper:
\begin{theorem} The product of degree $m$ Bernstein polynomials and its integral satisfies
\begin{eqnarray}
\prod_{j=1}^{m}B_{i_{j}}^{n_{j}}(t)=\frac{{M }}{{N \choose J}}B_{J}^{N}(t),
\end{eqnarray}
\begin{eqnarray}
\int_{0}^{1}\prod_{j=1}^{m}B_{i_{j}}^{n_{j}}(t)dt=\frac{M}{(N+1){N \choose J}},
\end{eqnarray}

\end{theorem}
where $M=\prod_{j=1}^{m}{n_{j} \choose i_{j}}$, $N=\sum_{j=1}^{m}n_{j}$ and $J=\sum_{j=1}^{m}i_{j}$.

\begin{theorem} Let the two polynomials $f(t)$ and $g(t)$ of degree $m$ and $n$ with coefficients $f_{i}^{m}$ and $g_{i}^{n}$ be as follows 

\begin{displaymath}
f(t)=\sum_{i=a}^{k}f_{i}^{m}B_{i}^{m}(t), \ \
g(t)=\sum_{i=0}^{n}g_{i}^{n}B_{i}^{n}(t),
\end{displaymath}
where $0\leq a<k\leq m$.
 Their product is a degree $m+n$ polynomial

\begin{eqnarray}
f(t)g(t)=\sum_{i=a}^{k+n} \left ( \sum_{j=max(a,i-n)}^{min(k,i)}\frac{{m \choose j}{n \choose i-j}}{{m+n \choose i}} f_{j}^{m}g_{i-j}^{n}\right)B_{i}^{m+n}(t).
\label{eq:4}
\end{eqnarray}
\end{theorem}
\begin{theorem} For each appropriate function $f(x)\in C_{[a,b]}$, there is a unique  polynomial $P_{n}(x)$  of degree at most $n$ to approximate $f(x)$, such that $\int_{a}^{b}(f(x)-P_{n}(x))^{2}dx$ is minimum \emph{(\citep {Isaacson})}.
\end{theorem}

\subsection{Statement of the approximation problem}
The problem of approximating  rational B\'{e}zier curves in equation (1) by constrained B\'{e}zier curves of arbitrary degree is that of finding  control points $\textbf{\textit{q}}_{0},\textbf{\textit{q}}_{1}, ...,  \textbf{\textit{q}}_{m}$, which define a  B\'{e}zier curve of degree $m$
\begin{eqnarray}
\textbf{Q}(t)=\sum \limits_{i=0}^{m}B_{i}^{m}(t)\textit{\textbf{q}}_{i},
\end{eqnarray}
such that the following two conditions are satisfied simultaneously:

1) B\'{e}zier curve $\textbf{Q}(t)$ has  the contact order $(k,h) \ (k,h=0,1,2)$ of  continuity at both endpoints of the rational B\'{e}zier curve $\textbf{\textit{P}}(t)$.

2) Given a weight function $\rho (t) >0$, a distance function $d_{\rho} (\textbf{\textit{P}},\textit{\textbf{Q}})$ between $\textbf{\textit{P}}(t)$ and $\textbf{\textit{Q}}(t)$ in the $L_{2}$ norm is minimized, that is
\begin{eqnarray}\nonumber
\begin{split}
   d_{\rho} (\textbf{\textit{P}},\textbf{\textit{Q}})&=\int_{0}^{1} \rho(t)\left(\textbf{\textit{P}}(t)-\textbf{\textit{Q}}(t)\right)^{2}dt. \\
&=\int_{0}^{1}\frac{\rho(t)}{\omega^{2}(t)}\left (\textbf{\textit{x}}(t)-\omega (t) \textbf{\textit{Q}}(t)\right)^{2}dt,
\label{eq:d}
\end{split}
\end{eqnarray}
and the control points $\textit{\textbf{q}}_{l} (l=k+1,k+2,...,m-h-1)$ of $\textbf{Q}(t)$ satisfy
\begin{eqnarray}
\frac{\partial d_{\rho}(\textit{\textbf{P}},\textit{\textbf{Q}})}{\partial \textit{\textbf{q}}_{l}}=2\int_{0}^{1}\left(\textbf{\textit{x}}(t)-\omega(t)\textbf{\textit{Q}}(t)\right)B_{l}^{m}(t)\frac{\rho(t)}{\omega(t)}dt=0.
\end{eqnarray}
Clearly,  in this paper the  function $\omega (t)$  is a factor of function $\rho(t)$.

  In this paper,we shall study constrained polynomials approximation in the next section. For simplicity, assume that $\rho(t)=\omega(t)$ and $\rho(t)=\omega^{2}(t)$ respectively.

\section{Algorithm for constrained B\'{e}zier approximation}
The approximation algorithm is completed through a two-step process. First, we calculate the constrained control points of the approximation B\'{e}zier curve $\textbf{\textit{Q}}(t)$ by imposing the contact order in the condition 1). Next, we calculate other unconstrained control points by using the weighted least-squares method, in order to satisfy condition 2).
\subsection{Constrained conditions of the approximation curve}
According to  equation (\ref{eq:energy}), we have
\begin{displaymath}
\textbf{\textit{x}}(t)=\textbf{\textit{P}}(t)\omega(t).
\end{displaymath}
Taking the derivative of the above equation $r (=0,1,2)$ times, we obtain
\begin{displaymath}
\textbf{\textit{x}}^{(r)}(t)=\sum \limits_{j=0}^{r}{r \choose j}\omega^{(j)}(t)\textbf{\textit{P}}^{(r-j)}(t).
\end{displaymath}
Thus, we have the zero, first and second derivatives of $\textbf{\textit{P}}(t)$ at two endpoints $(t=0,1)$,  satisfying the following,  respectively
\begin{displaymath}
\begin{split}
\textit{\textbf{P}}(0)&=\textit{\textbf{p}}_{0}, \ \textbf{\textit{P}}(1)=\textbf{\textit{p}}_{n}, \\
\textit{\textbf{P}}\ '(0)&=n\omega_{1}(\textbf{\textit{p}}_{1}-\textbf{\textit{p}}_{0}), \ \textit{\textbf{P}}\ '(1)=n\omega_{n-1}(\textbf{\textit{p}}_{n}-\textbf{\textit{p}}_{n-1}),
\\
\textit{\textbf{P}}\ ''(0)&=n\left \{\omega_{2}(n-1)\textbf{\textit{p}}_{2}+2\omega_{1}(1-n\omega_{1})\textbf{\textit{p}}_{1}+[2\omega_{1}(n\omega_{1}-1)+\omega_{2}(1-n)]\textbf{\textit{p}}_{0}\right\},
   \\
  \textit{\textbf{P}}\ ''(1)&=n  \{[2\omega_{n-1}(n\omega_{n-1}-1)+\omega_{n-2}(1-n)]\textbf{\textit{p}}_{n}+
  2\omega_{n-1}(1-n\omega_{n-1})\textbf{\textit{p}}_{n-1} \\
  &+(n-1)\omega_{n-1}\textbf{\textit{p}}_{n-2}  \}.
\end{split}
\end{displaymath}

Since the zero, first and second derivatives of $\textbf{\textit{Q}}(t)$ at two endpoints $(t=0,1)$ can be described as
\begin{displaymath}
\begin{split}
\textit{\textbf{Q}}(0)&=\textit{\textbf{q}}_{0}, \ \textbf{\textit{Q}}(1)=\textbf{\textit{q}}_{m}, \\
\textit{\textbf{Q}}\ '(0)&=m(\textbf{\textit{q}}_{1}-\textbf{\textit{q}}_{0}), \ \textit{\textbf{Q}}\ '(1)=m(\textbf{\textit{q}}_{m}-\textbf{\textit{q}}_{m-1}),
\\
\textit{\textbf{Q}}\ ''(0)&=m(m-1)(\textbf{\textit{q}}_{2}-2\textbf{\textit{q}}_{1}+\textbf{\textit{q}}_{0}),
 \\
  \textit{\textbf{Q}}\ ''(1)&=m(m-1)(\textbf{\textit{q}}_{m}-2\textbf{\textit{q}}_{m-1}+\textbf{\textit{q}}_{m-2}),
\end{split}
\end{displaymath}
matching the function value and derivatives up to the second order at both endpoints of $\textbf{\textit{P}}(t)$ and $\textbf{\textit{Q}}(t)$, the constrained control points of $\textbf{\textit{Q}}(t)$ are
  \begin{displaymath}
\begin{aligned}
\textbf{\textit{q}}_{0}&=\textbf{\textit{p}}_{0}, \textbf{\textit{q}}_{m}=\textit{\textbf{p}}_{n}, \\
\textrm{if}  \  m\geq 2, \ \ \ \ \textbf{\textit{q}}_{1}&=\frac{n}{m}\omega_{1}\textbf{\textit{p}}_{1}+\left(1-\frac{n}{m}\omega_{1}\right)\textbf{\textit{p}}_{0}, \\ \textrm{if} \  m \geq 3 , \textbf{\textit{q}}_{m-1}&=\left(1-\frac{n}{m}\omega_{n-1}\right)\textbf{\textit{p}}_{n}+\frac{n}{m}\omega_{n-1}\textbf{\textit{p}}_{n-1}, \\
\textrm{if} \  m \geq 4, \  \ \ \ \textbf{\textit{q}}_{2}&=\frac{1}{m(m-1)}\{n(n-1)\omega_{2}\textbf{\textit{p}}_{2}+2n\omega_{1}(m-n\omega_{1})\textbf{\textit{p}}_{1}
  \\
  &+[m(m-1)+2n\omega_{1}(n\omega_{1}-m)+n\omega_{2}(1-n)]\textbf{\textit{p}}_{0}  \}, \\
\textrm{if} \  m \geq 5 ,\textbf{\textit{q}}_{m-2}&=\frac{1}{m(m-1)}\{[2n\omega_{n-1}(n\omega_{n-1}-m)+n\omega_{n-2}(1-n)+m(m-1)]\textbf{\textit{p}}_{n} \\
&+{2n\omega_{n-1}}(m-n\omega_{n-1})\textbf{\textit{p}}_{n-1}+{n(n-1)}\omega_{n-2}\textbf{\textit{p}}_{n-2} \}.
\end{aligned}
\end{displaymath}
Accordingly, \textbf{\textit{Q}}(t) can be rewritten as
\begin{eqnarray}
\textbf{\textit{Q}}(t)=\sum \limits_{i=0}^{k}B_{i}^{m}(t)\textbf{\textit{q}}_{i}+\sum \limits_{i=k+1}^{m-h-1}B_{i}^{m}(t)\textbf{\textit{q}}_{i}+\sum \limits_{i=m-h}^{m}B_{i}^{m}(t)\textit{\textbf{q}}_{i},
\end{eqnarray}
where $\textbf{\textit{q}}_{i}(i=k+1,k+2, ..., m-h-1)$ are the unknown control points and
\begin{eqnarray*}
k=
\begin{cases}
0  & m=1, \\
1  & m=2,3, \\
2  & m\geq 4,
\end{cases}
\ \
h=
\begin{cases}
0  & m=1,2, \\
1  & m=3,4, \\
2  & m\geq 5 .
\end{cases}
\end{eqnarray*}
\subsection{Unconstrained control points of the approximation curve}
\subsubsection {$\rho (t)=\omega(t)$}
  When $\rho (t)=\omega(t)$, based on equation (6) we have
\begin{eqnarray}
\int_{0}^{1} \textbf{\textit{Q}}(t)\omega(t)B_{l}^{m}(t)dt=\int_{0}^{1}\textbf{\textit{x}}(t)B_{l}^{m}(t)dt.
\label{eq:dri}
\end{eqnarray}
Substituting (7) into (8), we get
\begin{eqnarray}
&&\sum \limits_{i=0}^{k}\int_{0}^{1}\omega(t)B_{i}^{m}(t)B_{l}^{m}(t)\textbf{\textit{q}}_{i}dt
+\sum \limits_{i=k+1}^{m-h-1}\int_{0}^{1}\omega(t)B_{i}^{m}(t)B_{l}^{m}(t)\textbf{\textit{q}}_{i}dt  \\
&&+\sum \limits_{i=m-h}^{m}\int_{0}^{1}\omega(t)B_{i}^{m}(t)B_{l}^{m}(t)\textbf{\textit{q}}_{i}dt \nonumber
=\int_{0}^{1}\textbf{\textit{x}}(t)B_{l}^{m}(t)dt.
\end{eqnarray}
By deriving  equation (9) based on  \textbf{theorems} 1 and 2, we obtain

\begin{eqnarray}
&&\sum \limits_{i=0}^{n+k}\left(\sum \limits_{j=max(0,i-n)}^{min(k,i)}\frac{{m \choose j}{n \choose i-j}}{{2m+n \choose i+l}}\omega_{i-j}\textbf{\textit{q}}_{j}\right)
+\sum \limits_{i=k+1}^{n+m-h-1}\left(\sum \limits_{j=max(k+1,i-n)}^{min(m-h-1,i)}\frac{{m \choose j}{n \choose i-j}}
{{2m+n \choose i+l}}\omega_{i-j}\textbf{\textit{q}}_{j}\right) \nonumber \\
&&+\sum \limits_{i=m-h}^{n+m}\left(\sum \limits_{j=max(m-h,i-n)}^{min(m,i)}\frac{{m \choose j}{n \choose i-j}}{{2m+n \choose i+l}}\omega_{i-j}\textbf{\textit{q}}_{j}\right) \nonumber \\
&&=\frac{2m+n+1}{m+n+1}\sum \limits_{i=0}^{n}\left(\frac{{n \choose i}}{{n+m \choose i+l}}\omega_{i}\textbf{\textit{p}}_{i}\right), \nonumber \\
&&(l=k+1,k+2,...,m-h-1).
\end{eqnarray}

To easily evaluate the unknown points $\textit{\textbf{q}}_{i} (i=k+1,k+2,...,m-h-1)$, we rearrange  equation (10) to obtain
\begin{eqnarray}
&&\sum \limits_{i=k+1}^{m-h-1}\sum \limits_{j=0}^{n}\frac{{n \choose j}{m \choose i}}{{2m+n \choose i+j+l}}\omega_{j}\textbf{\textit{q}}_{i}=
\frac{2m+n+1}{m+n+1}\sum \limits_{i=0}^{n}\frac{{n \choose i}}{{m+n \choose i+l}}\omega_{i}\textbf{\textit{p}}_{i} \nonumber \\
&&-\sum \limits_{i=0}^{k}\sum \limits_{j=0}^{n}\frac{{n \choose j}{m \choose i}}{{2m+n \choose i+j+l}}\omega_{j}\textbf{\textit{q}}_{i}-\sum \limits_{i=m-h}^{m}\sum \limits_{j=0}^{n}\frac{{n \choose j}{m \choose i}}{{2m+n \choose i+j+l}}\omega_{j}\textbf{\textit{q}}_{i},
\end{eqnarray}

$(l=k+1,k+2,...,m-h-1).$
\subsubsection{$\rho (t)=\omega(t)^{2}$}
 When $\rho (t)=\omega(t)^{2}$, one has
\begin{eqnarray}
\int_{0}^{1} \textbf{\textit{Q}}(t)\omega^{2}(t)B_{l}^{m}(t)dt=\int_{0}^{1}\textbf{\textit{x}}(t)\omega(t)B_{l}^{m}(t)dt.
\label{eq:dri}
\end{eqnarray}
Substituting (7) into (\ref{eq:dri}), it yields
\begin{eqnarray}
&&\sum \limits_{i=0}^{k}\int_{0}^{1}\omega^{2}(t)B_{i}^{m}(t)B_{l}^{m}(t)\textbf{\textit{q}}_{i}dt
+\sum \limits_{i=k+1}^{m-h-1}\int_{0}^{1}\omega^{2}(t)B_{i}^{m}(t)B_{l}^{m}(t)\textbf{\textit{q}}_{i}dt  \\
&&+\sum \limits_{i=m-h}^{m}\int_{0}^{1}\omega^{2}(t)B_{i}^{m}(t)B_{l}^{m}(t)\textbf{\textit{q}}_{i}dt \nonumber
=\int_{0}^{1}\textbf{\textit{x}}(t)\omega(t)B_{l}^{m}(t)dt.
\end{eqnarray}
By deriving  equation (13) based on  \textbf{theorems} 1 and 2, we obtain

\begin{eqnarray}
&&\sum \limits_{i=0}^{2n+k}\left(\sum \limits_{j=max(0,i-2n)}^{min(k,i)}\frac{{m \choose j}{2n \choose i-j}}{{2m+2n \choose i+l}}W_{i-j}\textbf{\textit{q}}_{j}\right)
+\sum \limits_{i=k+1}^{2n+m-h-1}\left(\sum \limits_{j=max(k+1,i-2n)}^{min(m-h-1,i)}\frac{{m \choose j}{2n \choose i-j}}
{{2m+2n \choose i+l}}W_{i-j}\textbf{\textit{q}}_{j}\right) \nonumber \\
&&+\sum \limits_{i=m-h}^{2n+m}\left(\sum \limits_{j=max(m-h,i-2n)}^{min(m,i)}\frac{{m \choose j}{2n \choose i-j}}{{2m+2n \choose i+l}}W_{i-j}\textbf{\textit{q}}_{j}\right) \nonumber \\
&&=\frac{2m+2n+1}{m+2n+1}\sum \limits_{i=0}^{2n}\left(\sum \limits_{j=max(0,i-n)}^{min(n,i)}\frac{{n \choose j}{n \choose i-j}}{{2n+m \choose i+l}}\omega_{i-j}\omega_{j}\textbf{\textit{p}}_{j}\right), \nonumber \\
&&(l=k+1,k+2,...,m-h-1),
\end{eqnarray}
where
\begin{displaymath}
W_{i}=\sum \limits_{j=max(0,i-n)}^{min(n,i)}\frac{{n \choose j}{n \choose i-j}}{{2n \choose i}}\omega_{i-j}\omega_{j}.
\end{displaymath}
We rearrange  equation (14) to obtain
\begin{eqnarray}
&&\sum \limits_{i=k+1}^{m-h-1}\sum \limits_{j=0}^{2n}\frac{{2n \choose j}{m \choose i}}{{2n+2m \choose i+j+l}}W_{j}\textbf{\textit{q}}_{i}=
\frac{2m+2n+1}{2m+n+1}\sum \limits_{i=0}^{n}\sum \limits_{j=0}^{n}\frac{{n \choose i}{n \choose j}}{{2n+m \choose i+j+l}}\omega_{i}\omega_{j}\textbf{\textit{p}}_{i} \nonumber \\
&&-\sum \limits_{i=0}^{k}\sum \limits_{j=0}^{2n}\frac{{2n \choose j}{m \choose i}}{{2n+2m \choose i+j+l}}W_{j}\textbf{\textit{q}}_{i}-\sum \limits_{i=m-h}^{m}\sum \limits_{j=0}^{2n}\frac{{2n \choose j}{m \choose i}}{{2n+2m \choose i+j+l}}W_{j}\textbf{\textit{q}}_{i},
\end{eqnarray}

$(l=k+1,k+2,...,m-h-1).$

Finally,  by \textbf{theorem} 3, we conclude that the solution of linear system (11) and (15) is unique. They can  be obtained by any methods introduced in
\citep {Isaacson}.

\section{Error estimation and implementation}

 Since
 \begin{eqnarray}\nonumber
 d_{\rho}(\textbf{\textit{P}},\textbf{\textit{Q}})&&=\left(\int_{0}^{1}\rho(t)\left|\textbf{\textit{P}}(t)-\textbf{\textit{Q}}(t)\right|^{p}dt\right)^{\frac{1}{p}}
 \\
 &&\leq\left(\int_{0}^{1}\rho(t)dt\right)^{\frac{1}{p}} \left( \max\limits_{0\leq t \leq 1} \left\|\textbf{\textit{P}}(t)-\textbf{\textit{Q}}(t)\right\| \right),
 \end{eqnarray}
 we use maximum distance $d_{max}$ to evaluate the approximation results for the convenience of estimation. More specifically, by dividing $[0, 1]$  into H subintervals of equal length $h = \frac{1}{H}$, it allows
\begin{displaymath}
d_{max}=\max\limits_{0\leq i \leq H}\|\textit{\textbf{P}}(t_{i})-\textbf{\textit{Q}}(t_{i})\|,
\end{displaymath}
where $t_{i} = ih \ (i = 0, 1, . . . , H)$ are evenly spaced values in the parameter domain $[0, 1]$.

 Another error for the approximation used here was given by Lu \citep{Lu}. That is
 $$d_{l_{1}}(\textbf{\textit{P}},\textbf{\textit{Q}})=\frac{h}{2}\sum_{i=0}^{H-1}(\|\textbf{\textit{P}}(t_{i})-\textbf{\textit{Q}}(t_{i})\|+\|\textbf{\textit{P}}(t_{i+1})-\textbf{\textit{Q}}(t_{i+1})\|).$$
 Obviously, $d_{max}>d_{l_{1}}$ by Equation (16). Finally, for more details about error analysis of weighted least-squares we refer the reader to {\cite{Watson}}.
 \\
 There is a disadvantage for approximating rational B\'{e}zier curve by  B\'{e}zier curve: when an original rational B\'{e}zier curve is convex, but the resulting B\'{e}zier curve may be non-convex. To deal with this problem, we use higher degree B\'{e}zier curve to reduce the approximation error (see Example 3).

\newpage

\textbf{Example 1}.  A rational B\'{e}zier curve of degree 3 is defined by the control points (0, 0), (0.2, 1.5), (0.8, 1.5), (1,0) and the associated weights 1, 1.2, 1.5, 1. Table 1. lists  the errors  obtained by the weighted least-squares and  the Lee's method in \citep{Lee} for $n=4,5,6,7,8,$ respectively. From the comparison of the errors , our algorithm can mostly produce better approximation results than Lee's method.

\begin{table} [h]
\centering

\caption{ Errors approximation  to the rational B\'{e}zier curve of degree 3 }

\label{address}

\begin{tabular}{|>{\tiny}c|>{\tiny}c|>{\tiny}c|>{\tiny}c|>{\tiny}c|>{\tiny}c|>{\tiny}c|>{\tiny}c|>{\tiny}c|}\hline

$n$ & \multicolumn{4}{|>{\tiny}c|}{maximum error} &\multicolumn{4}{|>{\tiny}c|}{$L_{1}$-error}\\\cline{2-9}
    & $\rho=\omega(t)$ &  $\rho=\omega^{2}(t)$& \multicolumn{2}{|>{\tiny}c|}{Lee's  method}& $\rho=\omega(t)$ &$\rho=\omega^{2}(t)$ &\multicolumn{2}{|>{\tiny}c|}{Lee's  method}\\\cline{4-5}\cline{8-9}
    &               &              &  The 1st        &  The  2nd &               &             & The 1st &  The 2nd \\
    &               &              &   method (r=100)  &   method                 &               &             & method(r=100)                  &  method \\\cline{1-9}
4	&0.022907       &0.023911	   & 0.020312 &0.053134	&0.004237	    &0.004185	      &0.004679    &0.004295\\
5	&0.007390	    &0.007743	   & 0.006150 &0.024009	&0.001162	    &0.001149         &0.002205	   &0.001424\\
6   &0.002849       &0.002987      & 0.004098 &0.011423 &4.371940e-04   &4.321077e-04     &0.002058    &6.297696e-04     \\
7	&8.065401e-04	&8.469361e-04  & 0.003461 &0.004955	&1.057862e-04	&1.045051e-04	  &0.002046    &2.004255e-04\\
8	&2.870556e-04	&3.017434e-04  & 0.003400 &0.002312	&3.871330e-05	&3.827086e-05	 & 0.002030    &9.057967e-05\\\hline

\end{tabular}
\end{table}

\newpage

\textbf{Example 2} (Also Example 1 in \citep{Lu}). A rational B\'{e}zier curve of degree 7 is defined by the control points (0, 0), (0.5, 2), (1.5, 2), (2.5, .2), (3.5,.2), (4.5, 2), (5.5, 2), (6, 0) and the associated weights 1, 2, 1/3, 2, 2, 1/3, 2, 1. Table 2. lists  the  errors  obtained by the weighted least-squares and the Lu's method in \citep{Lu} for $n=5,6,...,18,$ respectively. By comparison for $L_{1}$ error , our algorithm can produce better approximation results than Lu's method as $n=7,8,...,18$. However, for maximum error, we find that $\rho(t)=\omega(t)$ method generated smaller error values than $\rho(t)=\omega^{2}(t)$ method, but it is reverse for $L_{1}$ error.

\begin{table} [h]

\centering

\caption{Errors approximation  to the rational B\'{e}zier curve of degree 7}

\label{address}

\begin{tabular}{|>{\tiny}c|>{\tiny}c|>{\tiny}c|>{\tiny}c|>{\tiny}c|>{\tiny}c|}\hline

$n$ & \multicolumn{2}{|>{\tiny}c|}{maximum error} &\multicolumn{3}{|>{\tiny}c|}{$L_{1}$-error}\\\cline{2-6}
    & $\rho=\omega(t)$ &  $\rho=\omega^{2}(t)$&$\rho=\omega(t)$ &$\rho=\omega^{2}(t)$ & Lu's method ($\lambda$=0.95)\citep{Lu}\\\hline
5	&0.095480	&0.097236	&0.049830	&0.049584	         &0.043439\\
6	&0.084855	&0.086855	&0.047010	&0.046837	         &0.040826\\
7	&0.042564	&0.043096	&0.011913	&0.011892	         &0.037923\\
8	&0.035610	&0.037288	&0.012254	&0.012140	         &0.035547\\
9	&0.035619	&0.037298	&0.012246	&0.012133	         &0.023335\\
10	&0.005224	&0.005420	&0.001600	&0.001578	         &0.016259\\
11	&0.005374	&0.005620	&0.001566	&0.001537	         &0.013806\\
12	&0.003154	&0.003335	&9.999698e-04	&9.814476e-04	         &0.012295\\
13	&0.002257	&0.002361	&3.797283e-04	&3.675078e-04	         &0.011220\\
14	&0.001222	&0.001307	&2.998330e-04	&2.944257e-04	         &0.009933\\
15	&0.001206	&0.001290	&2.802343e-04	&2.754957e-04	         &0.008897\\
16	&2.834610e-04	&3.006113e-04	&5.642958e-05	&5.503411e-05	         &0.006065\\
17	&2.687450e-04	&2.860145e-04	&5.593685e-05	&5.458823e-05	         &0.002309\\
18	&1.240139e-04	&1.326072e-04	&2.950490e-05	&2.879857e-05	         &0.001462\\\hline
\end{tabular}
\end{table}

\textbf{Example 3.} Consider a  4th-degree rational B\'{e}zier curve  $\textbf{\textit{P}}$(t) with control points (0, 0), (0.2, 1.5),(0.5, 1.0), (0.8, 1.5), (1, 0) and the associated weights 1, 0.06, 0.08, 0.05, 1.   The $n$th-degree, $n=6,7,...,23$,  B\'{e}zier curves $\textbf{\textit{Q}}(t)$ approximating $\textbf{\textit{P}}(t)$ with the contact order $(1, 1)$ of continuity at two endpoints.    Table 3. lists the error obtained by $\rho(t)=\omega(t)$ and $\rho(t)=\omega^{2}(t)$ for $n=6,7,...,25$, respectively. For both maximum error and $L_{1}$ error, we find that  $\rho(t)=\omega(t)$ method generated smaller error values than $\rho(t)=\omega^{2}(t)$ method. Fig 1. shows that the resulting $C^{(1,1)}$-continuity approximation curve of degree 6 is nonconvex at near zero. Fig 2. shows that the given curve and the resulting curve of degree 18  are almost the same.

 \begin{table} [h]

\centering

\caption{Errors approximation  to the rational B\'{e}zier curve of degree 4 }

\label{address}

\begin{tabular}{|>{\tiny}c|>{\tiny}c|>{\tiny}c|>{\tiny}c|>{\tiny}c|}\hline

$n$ & \multicolumn{2}{|>{\tiny}c|}{maximum error} &\multicolumn{2}{|>{\tiny}c|}{$L_{1}$-error}\\\cline{2-5}
    & $\rho=\omega(t)$ &  $\rho=\omega^{2}(t)$&$\rho=\omega(t)$ &$\rho=\omega^{2}(t)$\\\hline
 6	&0.044480318880273	&0.051688479122085	&0.029655280096461	&0.031476147607478\\
7	&0.023868501893973	&0.032255328068047	&0.015879504757276	&0.017049959327487\\
8	&0.015149900429568	&0.018966490675292	&0.009620768974042	&0.010354259252695\\
9	&0.007889349654775	&0.010828999111530	&0.005297861401833	&0.005727317319065\\
10	&0.005260808682043	&0.006853089108748	&0.003256996716366	&0.003527894117667\\
11	&0.002696794398629	&0.003727944822126	&0.001818736806966	&0.001974652524091\\
12	&0.001847253697888	&0.002461821331414	&0.001126845692919	&0.001225098891989\\
13	&9.389817537108827e-04	&0.001303331966052	&6.345106044598648e-04	&6.908760850927994e-04\\
14	&6.534387865426517e-04	&8.830834863300054e-04	&3.949335561005876e-04	&4.303886131207782e-04\\
15	&3.306456957305267e-04	&4.601136641000896e-04	&2.236082280218562e-04	&2.439441476699068e-04\\
16	&2.323377554374534e-04	&3.169111231735876e-04	&1.395669996426054e-04	&1.523518332243233e-04\\
17	&1.172969349325125e-04	&1.635210281226485e-04	&7.933336636596293e-05	&8.667628246221323e-05\\
18	&8.293326487723877e-05	&1.138220713997069e-04	&4.960554098475324e-05	&5.421467244232657e-05\\
19	&4.182337819389267e-05	&5.838233410876038e-05	&2.828053820034141e-05	&3.093173320810643e-05\\
20	&2.968719111016691e-05	&4.092516147539204e-05	&1.770286264378748e-05	&1.936610741449082e-05\\
21	&1.496952671851878e-05	&2.091183351643694e-05	&1.011718721185016e-05	&1.106881100291165e-05\\
22	&1.064108236597488e-05	&1.472951815919057e-05	&6.336764567038676e-06	&6.938302189360695e-06\\
23	&5.374309134732559e-06	&7.516576177590380e-06	&3.631203840035559e-06	&3.974634105914670e-06\\\hline

\end{tabular}
\end{table}

\begin{figure}[tbp]
\begin{center}
\includegraphics[width=0.4\columnwidth]{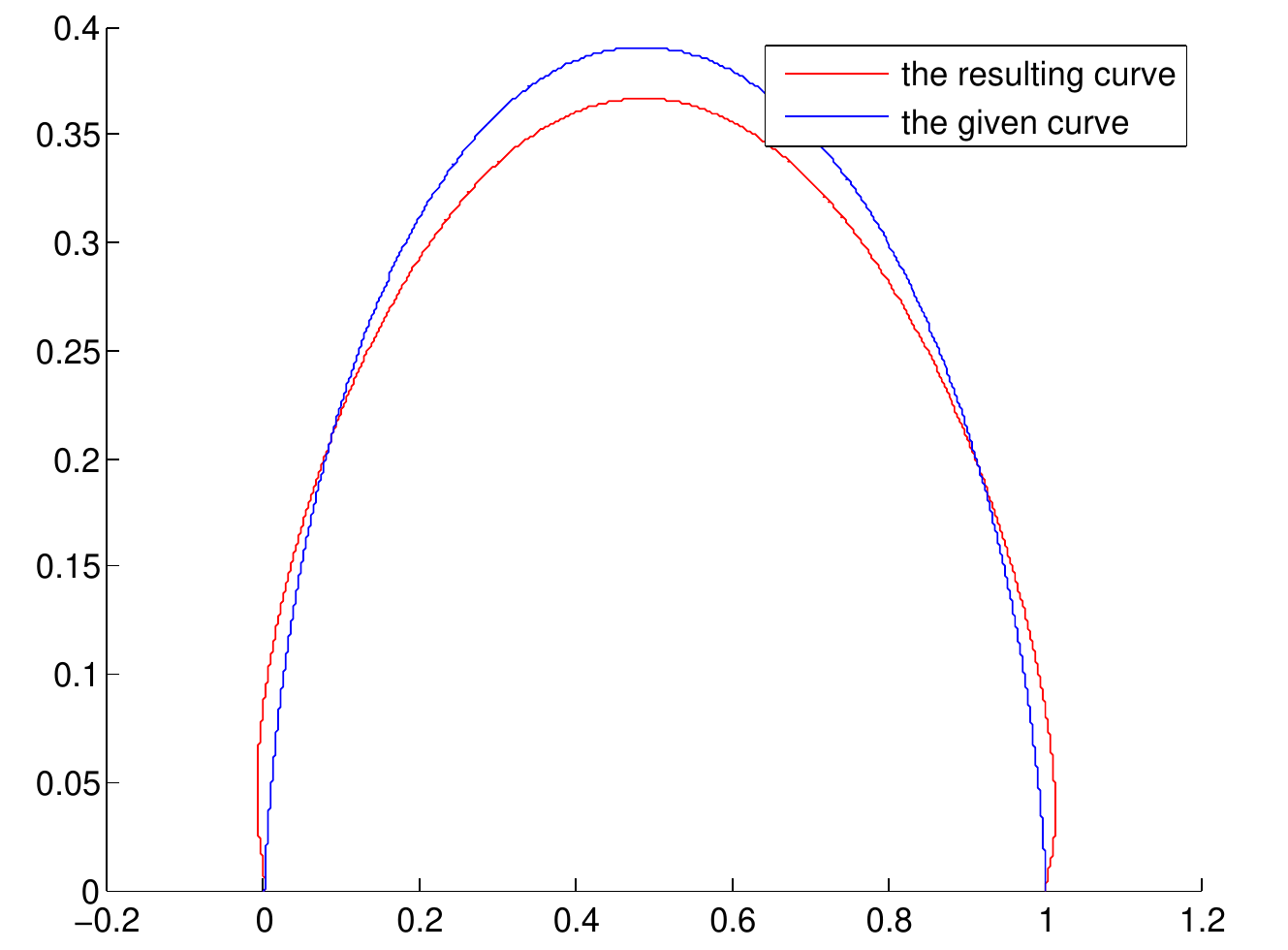} \ \includegraphics[width=0.4\columnwidth]{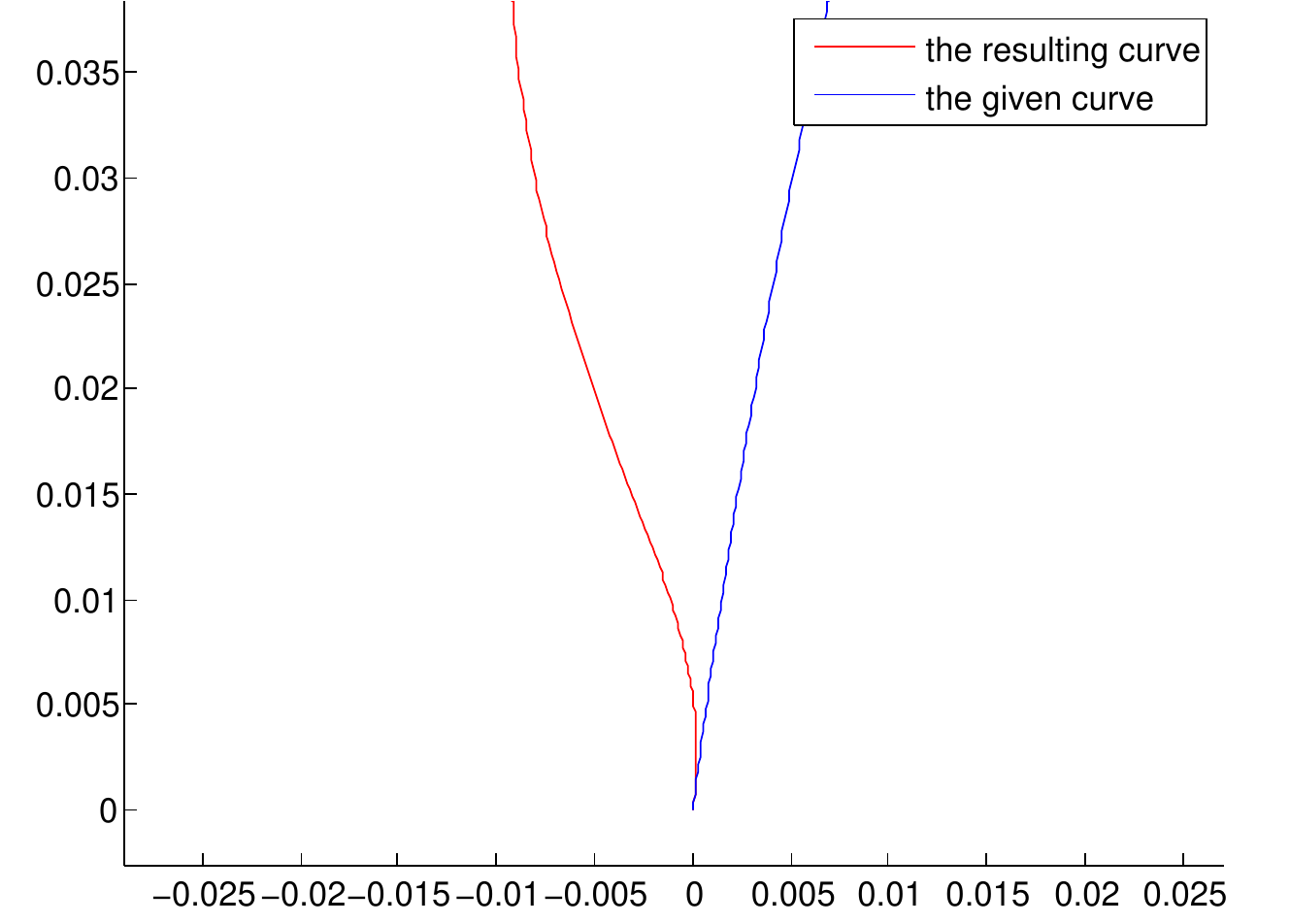}\\
(a) \hspace{5cm} (b)
\caption{ (a) The rational B\'{e}zier curve of degree 4 and the resulting $C^{(1,1)}$ approximation curve of degree 6 as $\rho(t)=\omega(t)$, (b) The photomicrograph of (a) at near zero.}
\label{fig1}
\end{center}
\end{figure}

\begin{figure}[tbp]
\begin{center}
\includegraphics[width=0.4\columnwidth]{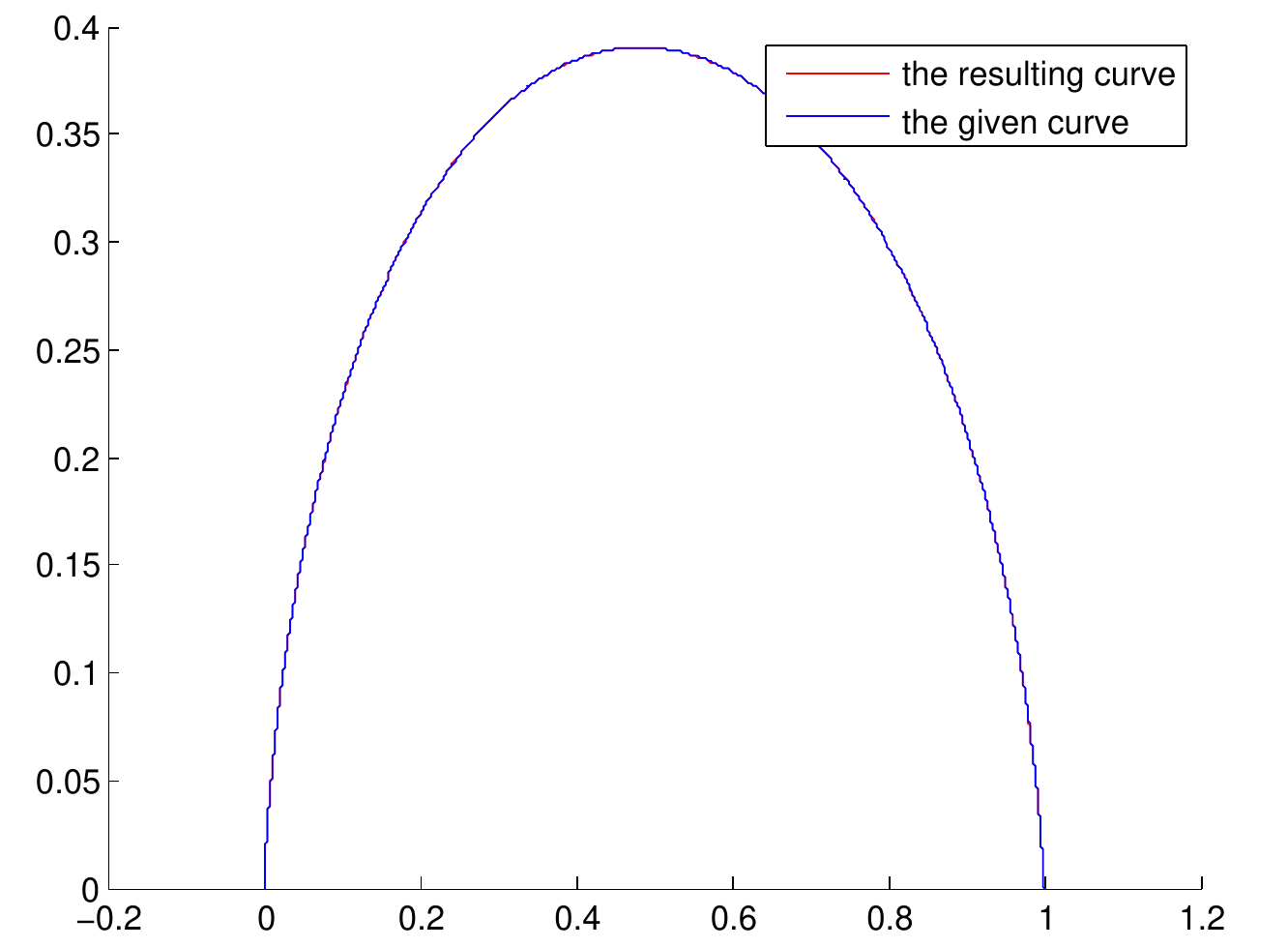} \includegraphics[width=0.4\columnwidth]{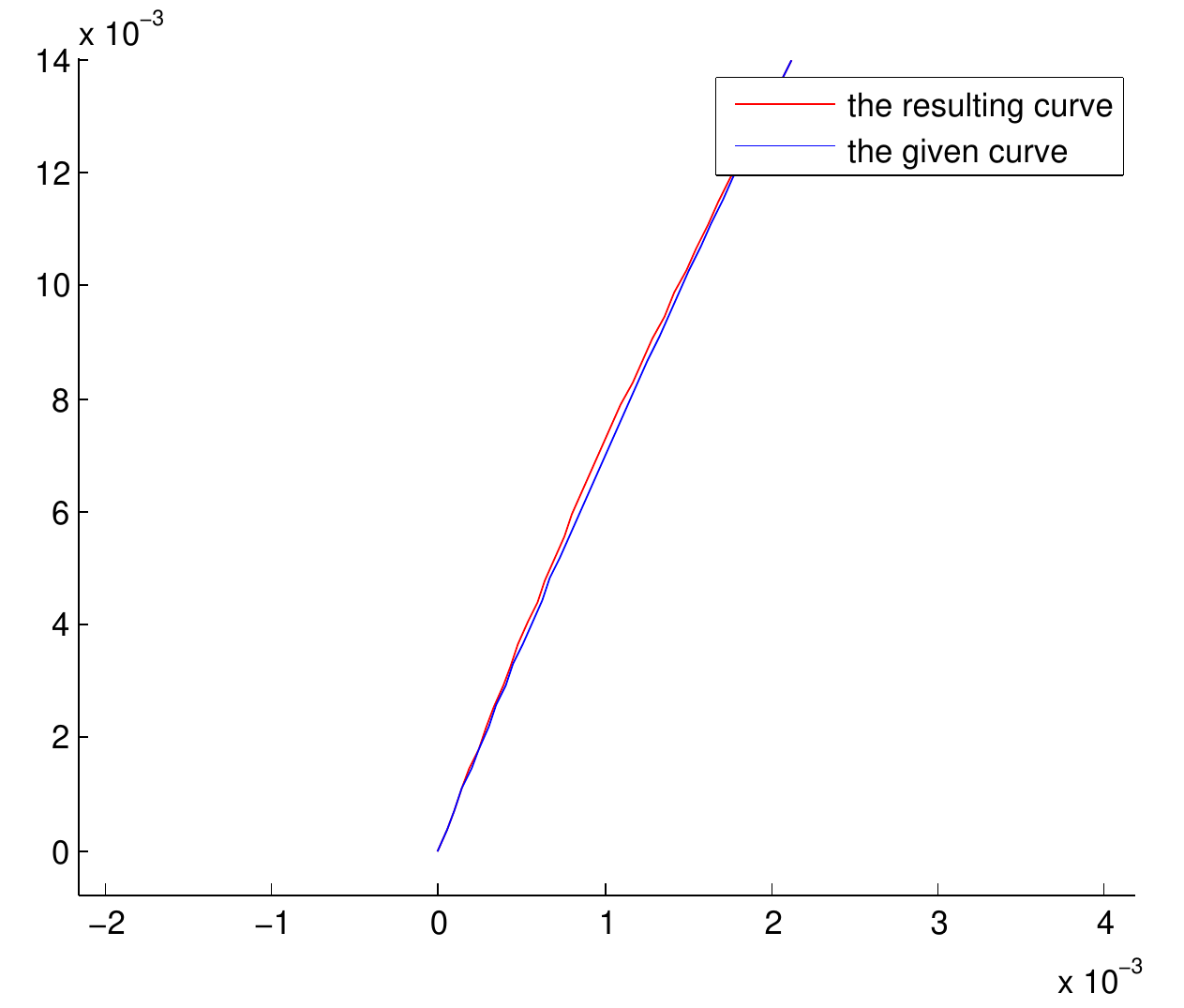}\\
(a) \hspace{5cm} (b)
\caption{ (a) The rational B\'{e}zier curve of degree 4 and the resulting $C^{(1,1)}$ approximation curve of degree 18 as $\rho(t)=\omega(t)$, (b) The photomicrograph of (a) at near zero.}
\label{fig1}
\end{center}
\end{figure}

\newpage
\textbf{Example 4.} (Also Example 1 in \citep{Hu13}).
The given curve is a rational B\'{e}zier curve of degree $8$ with the control points $(0, 0)$, $(0, 2)$, $(2, 10)$, $(4, 6)$, $(6, 6)$,$(11, 16)$, $(8, 1)$, $(9, 1)$, $(10, 0)$ and the associated weights $1$, $2$, $ 3$, $ 9$, $12$, $20$, $30$, $4$, $ 1$. We find a 5th-degree B\'{e}zier curve with different continuity to approximate the given curve, see Table 4. Here, the Hausdorff distance between curves $\textit{\textbf{P}}(t)$ and $\textit{\textbf{Q}}(t)$ is used. As shown in Table 4. our method is better than Hu's method for $C^{(1,1)}$-continuity when $\rho(t)=\omega^{2}(t)$.
\begin{table} [h]

\centering

\caption{Errors approximation  to the rational B\'{e}zier curve of degree 8 }

\label{address}

\begin{tabular}{|>{\tiny}c|>{\tiny}c|>{\tiny}c|>{\tiny}c|}\hline

$C^{(\alpha,\beta)}-\textmd{continuity}$ & \multicolumn{3}{|>{\tiny}c|}{Hausdorff distance errors}\\\cline{2-4}
    &   Hu's method &$\rho=\omega(t)$ &$\rho=\omega^{2}(t)$\\\hline
 $C^{(0,0)}$		&0.245371($\lambda$=1.046971)	& 0.376676194034764	& 0.487647609642435\\
$C^{(1,1)}$		&0.560612 ($\lambda$=0.713693)	&0.576850674814240	&0.519236375672234\\\hline

\end{tabular}
\end{table}
\newpage

\section{Conclusion}
In this paper, we have proposed an algorithm for approximating rational B\'{e}zier curves by weighted least-squares. From examples, we can find that approximation methods for $\rho(t)=\omega(t)$ and $\rho(t)=\omega^{2}(t)$ have their advantages and disadvantages respectively. To obtain convex-reserving approximation and reduce the approximation error, higher degree B\'{e}zier curve is used in this paper. These two problems are still our future works. Especially, orthogonal polynomials are applied to solve  problems .

\bigskip
\noindent {\bf Acknowledgments.}
The work is supported  by the  Natural Science Basic Research Plan in Shaanxi Province of China (No.2013JM1004),the  Fundamental Research Funds for the Central Universities (No.GK201102025), the National Natural Science Foundation of China (No.11101253) and the Starting Research Fund from the Shaanxi Normal University (No.999501).



\end{document}